\documentclass[12pt]{article}

\usepackage[margin=1in]{geometry}
\usepackage{graphicx}
\usepackage[nosort]{cite}
\usepackage{amsmath,amsthm,amsfonts}

\newtheorem{conjecture}{Conjecture}
\newtheorem*{mainconjecture}{Main Conjecture}
\newtheorem*{mainconjecturecontinued}{Main Conjecture, continued}

\numberwithin{equation}{section}

\newcommand{\Sing}{\operatorname{Sing}}
\newcommand{\ra}{\rangle}
\newcommand{\la}{\langle}
\newcommand{\CY}{Calabi--Yau}
\newcommand{\R}{\mathbb{R}}
\renewcommand{\P}{\mathbb{P}}

\newcommand{\td}{H}

\begin{document}

% \title[SLAG torus fibrations]{Special 
\title{Special
Lagrangian torus fibrations of complete intersection
Calabi--Yau manifolds:\\ a geometric conjecture\\[6pt]
{\small To Shing-Tung Yau on his $65^{\text{th}}$ birthday}}

\author{David R. Morrison$^*$ and M. Ronen Plesser$^\dagger$}

\date{}

\maketitle

{\small

{}$^*$ Departments of Mathematics and Physics, U.C. Santa Barbara, Santa Barbara CA 93106

{}$^\dagger$ Center for Geometry and Theoretical Physics, Duke University,
Durham NC 27708

}

\bigskip

\begin{abstract}
For complete intersection Calabi--Yau manifolds in toric varieties, Gross
and Haase--Zharkov have given a conjectural
combinatorial description of the special
Lagrangian torus fibrations whose existence was predicted by Strominger,
Yau and Zaslow.  We present a geometric version of this construction,
generalizing an earlier conjecture of the first author.

\end{abstract}

Federer observed in his classic monograph \cite{MR0257325}
that a compact complex submanifold of a
K\"ahler manifold minimizes volume in its homology class, due to a local
property of the K\"ahler form which follows from the Wirtinger inequality.
Harvey and Lawson \cite{MR666108}
generalized this to the notion of {\em calibrated
submanifold}\/ in 1982, and found several other classes of submanifolds
with a homological volume-minimizing property, 
including  the {\em special Lagrangian submanifolds
of a Ricci-flat K\"ahler manifold.}  Thanks to Yau's proof \cite{Yau,MR480350}
of the Calabi conjecture \cite{Calabi}, we have a rich supply of 
such Ricci-flat manifolds.

A Calabi--Yau manifold $Z$ has both a Ricci-flat K\"ahler form 
$\omega_Z$ and a 
non-vanishing holomorphic form of top degree $\Omega_Z$; the special
Lagrangian condition on a submanifold $L\subset Z$ 
(whose real dimension
is half that of $Z$)
simply says that both $\omega_Z$ and 
$\operatorname{Re}(\Omega_Z)$ vanish when restricted to $L$.
In spite of the simplicity of the definition,
the structure of special Lagrangian submanifolds is
still largely unknown.  However, given a special Lagrangian submanifold $L$
in a fixed Calabi--Yau manifold $Z$, the local 
deformation theory \cite{MR1664890} and global deformation theory
\cite{dg-ga/9711002} of $L$ within $Z$ are known, and in particular
it is known that the deformation space is a real manifold of 
dimension $b_1(L)$.
Thus, if $Z$ has complex dimension $n$ and $L$ has the topology of
a real $n$-torus $T^n$ and if moreover the nearby deformations of $L'$
are disjoint from $L$, then the resulting family of tori will determine
a fibration of an open subset of $Z$.

Nearly twenty years ago, Strominger, Yau, and Zaslow \cite{syz}
conjectured a relationship between the phenomenon of mirror
symmetry which had been discovered in the physics community
\cite{Aspinwall:1990xe,CLS,Greene:1990ud}
and fibrations of Calabi--Yau manifolds by special Lagrangian tori
(with singular fibers allowed).
According to the conjecture, special Lagrangian torus fibrations
should exist for Calabi--Yau manifolds which have mirror partners, and
the special Lagrangian torus fibrations on a mirror pair of 
Calabi--Yau manifolds should be dual to each other fiber by
fiber (at least once certain corrections to the geometry have
been made which are associated to holomorphic disks whose boundary
lies on a special Lagrangian torus).
The original physics argument is expected to apply when the moduli of the
Calabi--Yau metric are near the boundary of the moduli space.
More precisely, the complex
structure should be near a degeneration with ``maximally unipotent
monodromy'' \cite{mirrorguide,delignelimit,compact} and the K\"ahler class should be deep within the K\"ahler
cone.  
By restricting attention to small neighborhoods of the boundary,
and focussing on properties of the base of the fibration,
Gross and Siebert were
led to a beautiful reformulation of the conjecture
 as a problem in algebraic geometry \cite{GrossSiebert:log}, and much progress
has been made on that 
reformulation \cite{GrossSiebert:logI,GS:mslogII,GS:affine-complex}.\footnote{For a recent
review, see \cite{MR3204345}.}

\section{Geometric conjectures}

Returning to the original version of the problem,
even almost twenty years later we still lack the analytic tools to directly
analyze special Lagrangian submanifolds of a compact Calabi--Yau manifold.
However, building on some local analysis of Joyce 
\cite{Joyce:SYZ,Joyce:U1I,Joyce:U1II,Joyce:U1III} as well as early
work on the problem by Zharkov \cite{Zharkov:torus}, Gross
\cite{Gross:slagI,Gross:slagII,Gross:Tmir} and Ruan
\cite{Ruan:quinticI,Ruan:quinticII,Ruan:quinticIII,Ruan:different,Ruan:lagmir,Ruan:Newton,Ruan:hypersurfacesI,Ruan:hypersurfacesII,Ruan:hypersurfacesIII},
the first author
formulated in \cite{susyT3} a series of conjectures about
the structure of these fibrations.  Following \cite{susyT3}, we only state the conjectures for Calabi--Yau threefolds.

The first conjecture, essentially due to Gross \cite{Gross:slagI,Gross:slagII,Gross:Tmir}
and Ruan \cite{Ruan:quinticI,Ruan:quinticII,Ruan:quinticIII},
concerns the combinatorial properties of the fibration.

\begin{conjecture}
\label{con:topology}
Let $\pi: Z\to B$ be a special Lagrangian $T^3$
fibration of a compact
Calabi--Yau threefold with respect to a Calabi--Yau metric
whose compatible complex structure is sufficiently close to a
boundary point with maximally unipotent monodromy,
and whose
K\"ahler class is sufficiently deep in the K\"ahler cone.  Then
\begin{itemize}

\item[i)] The discriminant locus of the fibration retracts onto
  a trivalent graph $D$, which we call the {\em combinatorial discriminant
locus.}

\item[ii)] For any loop around an edge of $D$, the monodromy on 
either $H^1\cong H_2$ or $H_1$ of the $3$-tori
is conjugate to
\begin{equation} \label{eq:mono3}
M=\left(\begin{array}{ccc} 1&0&1\\0&1&0\\0&0&1 \end{array} \right) \thinspace.
\end{equation}
In particular, both monodromy actions have a $2$-dimensional fixed plane.

\item[iii)] The vertices of $D$ come in two types: near a 
{\em positive vertex},
the three monodromy actions on $H^1\cong H_2$ near the vertex have 
fixed planes whose intersection is $1$-dimensional, 
while the three monodromy actions
on $H_1$ have a common $2$-dimensional fixed plane.
In an appropriate basis, the monodromy matrices on $H_1$ take the form
\begin{equation} \label{eq:posmat}
\left(\begin{array}{ccc} 1&0&1\\0&1&0\\0&0&1 \end{array} \right) ,
\quad
\left(\begin{array}{ccc} 1&0&0\\0&1&1\\0&0&1 \end{array} \right) ,
\quad
\left(\begin{array}{ccc} 1&0&-1\\0&1&-1\\0&0&\hphantom{-}1 \end{array} 
\right) .
\end{equation}
On the other hand, near a {\em negative vertex},
the three monodromy actions on $H^1\cong H_2$ near the vertex have a 
common $2$-dimensional fixed plane, while the three monodromy actions
on $H_1$ have fixed planes whose intersection is $1$-dimensional.
In an appropriate basis, the monodromy matrices on $H_1$ take the form
\begin{equation} \label{eq:negmat}
\left(\begin{array}{ccc} 1&0&1\\0&1&0\\0&0&1 \end{array} \right) ,
\quad
\left(\begin{array}{ccc} 1&1&0\\0&1&0\\0&0&1 \end{array} \right) ,
\quad
\left(\begin{array}{ccc} 1&-1&-1\\0&\hphantom{-}1&\hphantom{-}0\\0
&\hphantom{-}0&\hphantom{-}1 \end{array} \right) .
\end{equation}

\item[iv)] 
There are open neighborhoods $U_P\subset B$ of the vertices $P$ of $D$
such that each fiber of $\pi$ over a point not in $\bigcup_P U_P$ has
Euler characteristic $0$, and such that the Euler characteristic of
$\pi^{-1}(U_P)$ is either $1$ or $-1$, depending on whether
$P$ is a positive or negative vertex.

\end{itemize}
\end{conjecture}

The second conjecture, formulated in \cite{susyT3}, is more geometric
in nature.  It should be stressed that there is very little evidence
for this geometric conjecture at the present moment.  However, it appears
to be geometrically very natural, and assuming its truth has helped to
clarify a number of the tricky combinatorial issues associated with
mirror symmetry.

\begin{conjecture}
\label{conj:big}
Let $\pi: Z\to B$ be a special Lagrangian $T^3$
fibration of a compact
Calabi--Yau threefold with respect to a Calabi--Yau metric
whose compatible complex structure is sufficiently close to a
boundary point with maximally unipotent monodromy,
and whose
K\"ahler class is sufficiently deep in the K\"ahler cone.  Then
\begin{itemize}
\item[i)] The set $C\subset Z$ of singular points of fibers of $\pi$ is a
complex subvariety of $Z$ of complex dimension $1$.

\item[ii)] All singular points of $C$ are transverse triple points, 
locally of the
form $\{z_1z_2=z_1z_3=z_2z_3=0\}$ for local complex coordinates $z_1$,
$z_2$, $z_3$.

\item[iii)] For each connected component $C_\alpha$ of $C$, 
the image
$\pi(C_\alpha)$ has the topology of a disc with $g(C_\alpha)$ holes,
and the map
$\pi|_{C_\alpha}$ is generically $2$-to-$1$ onto its
image.  The image $\pi(C_\alpha)$ has a distinguished boundary
component which contains the images of all points in
$\Sing(C)\cap C_\alpha$.

\item[iv)] There is a  graph $D_\alpha\subset \pi(C_\alpha)$ 
to which $\pi(C_\alpha)$ retracts.
$D_\alpha$ has univalent
vertices on the boundary of $\pi(C_\alpha)$ at the images of
all points in $\Sing(C)\cap C_\alpha$, and has only
trivalent vertices in the interior.  The graph $D = \bigcup_\alpha D_\alpha$ is the combinatorial discriminant locus of Conjecture~\ref{con:topology}.

\item[v)] 
The map $\pi$ puts the singular points of $C$ in one-to-one 
correspondence with the positive vertices of $D$, which are
located at the intersections of the various $D_\alpha$.
The negative vertices of $D$
are the interior vertices of the various $D_\alpha$'s.

\end{itemize}
\end{conjecture}

\begin{figure}
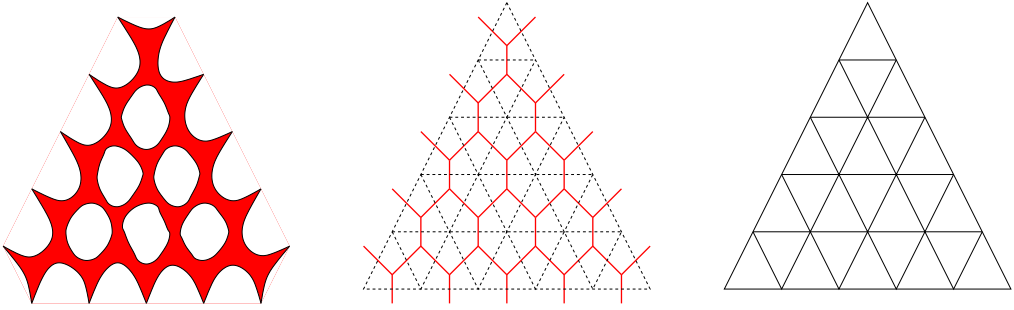

\begin{center}
\includegraphics[scale=0.3]{5triangle7.mps}
\qquad
\includegraphics[scale=0.3]{5triangle3b.mps}
\qquad
\includegraphics[scale=0.3]{5trianglea.mps}
\end{center}
\caption{Ingredients of the geometric conjecture for the quintic
threefold: the image of $C_{ij}$ under the moment map,
a graph to which it retracts, and a part of the corresponding
triangulation of the dual polytope for $\mathbb{CP}^2_{ij}$.}
\label{fig:quintic}
\end{figure}

The singular locus  $C$, its image $\pi(C)$ and the graph
$D$ were conjecturally described for Calabi--Yau hypersurfaces 
in toric varieties in \cite{susyT3}, motivated by a construction
of Zharkov \cite{Zharkov:torus}.  
We illustrate the conjecture in the case of quintic hypersurfaces in
$\mathbb{CP}^4$.  Let $x_1$, \dots, $x_5$ be the homogeneous coordinates
of $\mathbb{CP}^4$, and let $F(x_1,\dots,x_5)$ be the defining equation
of the hypersurface.  The components of $C$ are the quintic
plane curves
\[C_{ij}:=\{x_i=x_j=F=0\}\subset \mathbb{CP}^2_{ij}.\]
They meet along the sets
\[P_{ijk}:=\{x_i=x_j=x_k=F=0\},\]
each of which consists of five points.
The image $\pi(C_{ij})$ is isomorphic to the image of
$C_{ij}$ under
the {\em moment map}\/ $\mu_{ij}: \mathbb{CP}^2_{ij} \to \mathbb{R}^2$,
as illustrated on the left side of Figure~\ref{fig:quintic}.
The graph $D_{ij}$, illustrated in the middle of Figure~\ref{fig:quintic},
 is a tropical limit of 
$C_{ij}$ to which $\mu_{ij}(C_{ij})$ retracts;
it is determined by a triangulation of the dual polytope
of $\mathbb{CP}^2_{ij}$ (illustrated on the right
side of Figure~\ref{fig:quintic}).

The goal of this paper is to 
describe $C$, $\pi(C)$, and $D$
 for complete intersections in toric varieties.

\section{Nef partitions and singular limits}\label{sec:NP}

The combinatorial tool needed
 to describe a Calabi--Yau complete intersection in a toric
variety is a ``nef partition,'' first defined by Borisov \cite{borisov-ms}
and later expanded upon by 
Batyrev and Borisov \cite{MR1416334,Batyrev-Borisov-1994}.
Let $X=X_T$ be a compact toric variety of dimension $d$ determined 
by\footnote{The polytope determining $X$ is often denoted by $\Delta$
in the literature.  For the current construction, although we will indeed
encounter a polytope $\Delta\subset N_{\mathbb{R}}$, it is distinct from
$\nabla^\vee$ except in the hypersurface case.}
a triangulation $T$ of a polytope 
$\nabla^\vee \subset N_{\mathbb{R}}\cong \mathbb{R}^d$.  (The fan determining
$X_T$ consists of cones over simplices in $T$; alternatively, the vertices
of $T$ provide a construction of $X_T$ as a geometric quotient: see
\cite{Fulton:,cox}).

A vertex $e_a$ of $T$  determines a toric divisor  $\td_a$.  
A family of complete intersections $Z$ is determined by 
a {\em nef partition}\/ of $\nabla^\vee$ with
$r$ parts. This is a partition of the vertices $e_1,\ldots e_N$ 
into $r$ disjoint groups 
\begin{equation}
\{e_1,\ldots,e_N\} = \amalg_{i=1}^r E^{(i)}, \quad \text{where } E^{(i)} =
\{e_{i,1},\ldots, e_{i,d(i)}\}\ ,
\end{equation}
with the property that the associated line bundles 
\[ {\cal L}_i := {\cal O}_X\left(\sum_{\alpha=1}^{d_i} \td_{i,\alpha}\right)\]
are all nef.\footnote{This condition can be expressed
combinatorially \cite{borisov-ms}.}

We then have 
\begin{equation}
\nabla^\vee = {\rm Conv}\left(\Delta^{(1)}\cup\cdots\cup\Delta^{(r)}\right)\
,
\end{equation}
where 
\begin{equation}
\Delta^{(i)} = {\rm Conv}\left(\{0\}\cup E^{(i)}\right)\ .
\end{equation}
Note that the product of the nef line bundles
\[ \otimes_i{\cal   L}_i 
= {\cal O}_X\left(\td_1+\cdots + \td_N\right) = 
{\cal O}_X(-K_X) 
\]
is also nef, which implies that $\nabla^\vee$ is a relexive polyhedron.

Our Calabi--Yau $n$-fold
is given as the complete intersection 
\[Z=\{F_1=\dots=F_r=0\}\subset X,\]
 where
$F_i$ is a homogeneous polynomial of the same degree as 
$\prod_{\alpha=1}^{d(i)}x_{i,\alpha}$, i.e., it is a global section
of the line bundle ${\cal L}_i$.

We restrict attention here to the case of an
{\sl irreducible\/} nef partition.  A partition is irreducible if
there is no subset $\{i_1,\ldots ,i_k\}\subset\{1,\ldots ,r\}$ such
that $\Delta^{(i_)} + \cdots + \Delta^{(i_k)}$ contains 0 in its
interior.   Up to a possible refinement of the lattice any nef
partition is a direct sum of irreducibles \cite{HZ:T3III}.  The
\CY\ spaces associated to reducible partitions are thus discrete
quotients of products of complete intersections of lower dimension.

For any $\vec{\alpha}=(\alpha_1,\dots,\alpha_r)$ with $1\le \alpha_i\le d(i)$, we define
\[ X_{\vec{\alpha}} = \bigcap_{i=1}^r \td_{i,\alpha_i}.\]
Note that each $X_{\vec{\alpha}}$ is a 
(possibly empty) toric subvariety of $X$ of codimension 
$r$.

For each $X_{\vec{\alpha}}$ and each toric divisor $\td_{\ell,m}$
with $m\ne\alpha_\ell$, the intersection
$X_{\vec{\alpha}}\cap \td_{\ell,m}$ coincides with $X_{\vec{\alpha}}\cap X_{\vec{\beta}}$
where $\vec{\beta}$ is obtained from $\vec{\alpha}$ by 
\[ \beta_i = \begin{cases} \alpha_i &\text{ if $i\ne \ell$ } \\ m &\text{ if $i=\ell$}
\end{cases}.\]
It follows that if we remove
\[ X_{\vec{\alpha}} \cap \left(\bigcup_{\vec{\beta}\ne\vec{\alpha}} X_{\vec{\beta}}\right)\]
from $X_{\vec{\alpha}}$, the remainder is isomorphic to $(\mathbb{C}^*)^n$,
and hence is fibered by real $n$-tori.

We can form a natural family of Calabi--Yau varieties by taking the 
$i^{\text{th}}$ defining polynomial to be $tF_i+(1-t)
\prod_{\alpha=1}^{d(i)}x_{i,\alpha}$; as $t\to0$, this approaches
the ``large complex structure limit.''  By construction,
the limiting Calabi--Yau variety in this large complex structure limit
is
\[ Z_0 := \bigcup_{\vec{\alpha}} X_{\vec{\alpha}}.\]
Points on $Z_0$ which belong to only one $X_{\vec{\alpha}}$ are part of
the torus fibration.

As in the hypersurface case, and as in Zharkov's early 
work \cite{Zharkov:torus}, we expect
that the $T^n$ fibration on the large complex structure limit will deform
along with the complex structure to nearby Calabi--Yau manifolds.

We can now formulate the main new conjecture of this paper.

\bigskip

\begin{mainconjecture}
If $Z$ is close to the large complex
structure limit $Z_0$, then the $n$-torus fibration on $Z_0$
deforms to a fibration  $\pi:Z\to B$ by real $n$-tori.  When singular fibers
are included, the parameter space for the fibration is diffeomorphic
to $S^n$.

Moreover, the set of singular points of fibers of this
$T^n$-fibration is a complex subvariety of complex codimension two
in $Z$, described as follows.  For each part of the nef partition, there
is a subvariety $C^{(i)}$ of codimension two and a subvariety of
$P^{(i)}$ of codimension three, defined by
\[ C^{(i)} = \bigcup_{1\le \lambda_1<\lambda_2\le d(i)}
Z \cap \td_{i,\lambda_1} \cap \td_{i,\lambda_2} , \text{ and }\]
\[ P^{(i)} = \bigcup_{1\le\mu_1<\mu_2<\mu_3\le d(i)}
Z \cap \td_{i,\mu_1}\cap \td_{i,\mu_2} \cap \td_{i,\mu_3}.\]
The conjectured set of singular points of fibers of the $T^n$ fibration
is $C:=C^{(1)}\cup\cdots\cup C^{(r)}$.

\end{mainconjecture}

Note that for $Z$ generic, the components $C^{(i)}_{\lambda_1,\lambda_2}$
of $C^{(i)}$ are all nonsingular, and also that $P^{(i)}$ is the singular
locus of $C^{(i)}$; each component $P^{(i)}_{\mu_1,\mu_2,\mu_3}$
serves as the intersection of three of the components of $C^{(i)}$.
Note also that for $i\ne j$, $C^{(i)}\cap C^{(j)}$ has complex
codimension at least four in $Z$.

The singular set $C$ is a complex curve
in the case of Calabi--Yau threefolds.
We will describe $\pi(C)$ in the next section, and the combinatorial
discriminant $D$ to which it retracts in section~\ref{sec:CD}.

One easy consequence of this conjecture is a formula for the Euler
characteristic of $Z$ when $n\le3$.  We explain this in examples.

The case $n=1$ is not interesting because there is nothing of codimension
$2$, and the Calabi--Yau $1$-fold is always fibered by $T^1$'s over $S^1$,
with no singular fibers.

In the case $n=2$, $C$ is a collection of points.  
Each singular fiber contributes $1$ to the Euler characteristic, and the
total space is a K3 surface, so we expect precisely $24$ points in $C$.
Let us see how that works in various cases.

A complete intersection of degree $(d_1,\dots,d_k)$ in $\mathbb{P}^{k+2}$
(with each $d_k\ge2$) is a K3 surface if $\sum d_i=k+3$.  We divide the
homogeneous coordinates into $k$ groups with $d_j$ elements in the
$j^{\text{th}}$ group.  In the $j^{\text{th}}$ group, we should intersect
the K3 surface with one of $\binom{d_j}2$ pairs of hyperplanes, and
each such intersection will have $\prod d_i$ points.  Thus, the total
number of points is
\[ \left( \sum_j \binom{d_j}2 \right) \prod_i d_i.\]
Remarkably, this turns out to be $24$ in every case:

\begin{center}
\begin{tabular}{l|l|l|l}
  $(d_1,\dots,d_k)$ & $\sum_j \binom{d_j}2 $ & $\prod_id_i$ & \# points\\[3pt] \hline
$(4)$ & $\binom42=6$ & $4$ & $24$ \\[3pt]
$(3,2)$ & $\binom32 + \binom22 = 4$ & $3\cdot2=6$ & $24$ \\[3pt]
$(2,2,2)$ & $\binom22 + \binom22 + \binom22 = 3$ & $2\cdot2\cdot2=8$
& $24$ \\ 
\end{tabular}
\end{center}

Now we consider the case $n=3$.  Each $C^{(i)}_{\lambda_1,\lambda_2}$
is a nonsingular complex curve, and they meet three at a time along the sets
$P^{(i)}_{\mu_1,\mu_2,\mu_3}$.  Note that the various complex curves
$C^{(1)}, C^{(2)}, \cdots C^{(r)}$ are pairwise disjoint
(since their intersection has codimension $4$ on $Z$).

To compute the Euler characteristic, we
remove all of the intersection points
from each $C^{(i)}_{\lambda_1,\lambda_2}$, leaving us with a punctured
complex curve $\widetilde{C}^{(i)}_{\lambda_1,\lambda_2}$
with negative Euler characteristic.  This
retracts onto a graph consisting solely of negative vertices,
with the number of vertices being the absolute value of the
Euler characteritic of $\widetilde{C}^{(i)}_{\lambda_1,\lambda_2}$.
Since the fiber of each positive (resp.\ negative) vertex contributes 
Euler characteristic ${}{+}1$ (resp.\ ${}{-}1$) of the Calabi--Yau
threefold, the overall Euler characteristic is
\[ 
\sum \#(P^{(i)}_{\mu_1,\mu_2,\mu_3})
-\sum |\chi(\widetilde{C}^{(i)}_{\lambda_1,\lambda_2})|\ .
\]

For a complete intersection of degree $(4,2)$ in $\mathbb{P}^6$, we
divide the homogeneous coordinates into a group of $4$ and a group of $2$.
There are $\binom42=6$ complex curves from the first group, and 
$\binom22=1$ from 
the second group.  The complex curves are complete intersections of degree
$(4,2)$ in $\mathbb{P}^3$, and have genus $1+4\cdot2=9$ and
Euler characteristic $-2\cdot(4\cdot2)=-16$.
There are $\binom43=4$  toric subvarieties of codimension three
associated to the first group but none from the second group
(since $\binom23=0$); each meets the Calabi--Yau in $4\cdot2=8$
points.
In the first group, the complex curves are each punctured along $4-2$ toric
subvarieties (corresponding to the two remaining coordinates in the first group, out of $4$), so at a total of $(4-2)\cdot4\cdot2=16$ points.  
Thus, the Euler characteristic
of these punctured curves
is $-2\cdot(4\cdot2) 2-(4-2)\cdot(4\cdot2)=-4\cdot(4\cdot2)=-32$.
On the other hand, the Euler characteristic of the (unpunctured)
curve in the second group $-2\cdot(4\cdot2) - (2-2)\cdot(4\cdot2)
= -2\cdot(4\cdot2)=-16$, where we are thinking of it being punctured
along $2-2$ toric subvarieties.
The total Euler characteristic is thus
\[\left(\binom43 + \binom23\right)\cdot (4\cdot 2)
- \binom42 \cdot 4 \cdot (4\cdot2) - \binom22\cdot 2 \cdot(4\cdot 2)
= -176.\]

The analysis of other complete intersections is similar, and we indicate
each one by a parallel equation which indicates how the Euler characteristic
is calculated.  For the quintic hypersurface, the calculation is
\[ \binom53\cdot (5) - \binom52\cdot 5 \cdot (5) = -200.\]
For the complete  intersection of degree $(3,3)$ in $\mathbb{P}^6$, the
calculation is
\[ \left(\binom33 + \binom33\right) \cdot(3\cdot3)
-\binom32 \cdot 3\cdot(3\cdot3) - \binom32\cdot 3\cdot(3\cdot3)
=-144.\]
For the complete intersection of degree $(3,2,2)$ in $\mathbb{P}^7$, 
the calculation is
\[ \left(\binom33+\binom23+\binom23\right) \cdot(3\cdot2\cdot2)
-\binom32 \cdot 3 \cdot(3\cdot2\cdot2) 
-\binom22\cdot 2\cdot(3\cdot2\cdot2)
-\binom22\cdot 2\cdot(3\cdot2\cdot2)
=-144.\]
Finally, for the complete intersection of degree $(2,2,2,2)$ in $\mathbb{P}^8$, the calculation is
\[ \left(4 \times \binom23 \right) \cdot (2^4)
- 4 \times (\binom22\cdot 2 \cdot (2^4)) = -128.\]
Notice that in this last case, the first term is zero:
there are no toric subvarieties of
codimension three, and no positive vertices, in the calculation.

\section{The image under $\pi$, and a further limit}\label{sec:pi}

We expect each component of $\pi(C)$ to have
real codimension one in $S^n$, with  an 
amoeba-like structure analogous
to the one shown in the left side of Figure~\ref{fig:quintic}.
However, this is not easy to see directly.

\begin{figure}
\begin{center}
\includegraphics[scale=0.5]{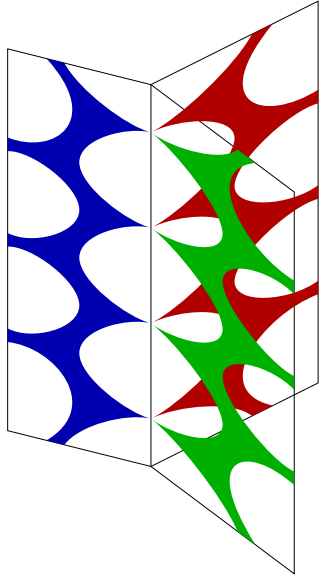}
\caption{Components of 
$\pi(C^{(i)})$ meeting three-at-a-time.} \label{fig:junction3}
\end{center}
\end{figure}

To describe our proposal for the structure
of $\pi(C)$, we need to consider a further
degeneration of the components
of the singular locus.  (This degeneration will also be useful
in section~\ref{sec:CD} in understanding the structure of the
combinatorial discriminant.)
Each of the components $C^{(i)}_{\lambda_1,\lambda_2}$ is itself 
a complete intersection in the toric variety $\td_{i,\lambda_1}\cap
\td_{i,\lambda_2}$, and we will take a further degeneration of
this complete intersection.  We do this by replacing each $F_j$, $j\ne i$
by the corresponding 
homogeneous monomial $\prod_{\alpha=1}^{d(j)}x_{j,\alpha}$.

For any $\vec{\alpha}=(\alpha_1,\dots,\alpha_r)$ with $1\le \alpha_i\le d(i)$, we define
\[ Y^{(i)}_{\vec{\alpha}} = \bigcap_{\substack{j=1\\j\ne i}}^r \td_{j,\alpha_j},\]
and note that each $Y^{(i)}_{\vec{\alpha}}$ is a 
(possibly empty) toric subvariety of $X$ of codimension 
$r{-}1$.  
Just as in the earlier analysis of $Z_0$, it is easy to see that
\[ \bigcap_{j\ne i} \left\{\prod_{\alpha=1}^{d(j)}x_{j,\alpha}=0\right\}
= \bigcup_{\vec{\alpha}} Y^{(i)}_{\vec{\alpha}}.\]
Our further degeneration thus replaces $C^{(i)}_{\lambda_1,\lambda_2}$
by
\[ \widehat{C}^{(i)}_{\lambda_1,\lambda_2} 
= \{F_i=0\} \cap H_{i,\lambda_1} \cap H_{i,\lambda_2} \cap 
\bigcup_{\vec{\alpha}} Y^{(i)}_{\vec{\alpha}} 
\]
whose components take the form
\[ \widehat{C}^{(i)}_{\lambda_1,\lambda_2,\vec{\alpha}}
= \{F_i=0\} \cap H_{i,\lambda_1} \cap H_{i,\lambda_2} \cap
Y^{(i)}_{\vec{\alpha}}.\]
These subvarieties meet along varieties
\[ \widehat{P}^{(i)}_{\mu_1,\mu_2,\mu_3,\vec{\alpha}}
= \{F_i=0\} \cap H_{i,\mu_1} \cap H_{i,\mu_2} \cap H_{i,\mu_3} \cap
Y^{(i)}_{\vec{\alpha}},\]
but they also meet along varieties
\[ \widehat{Q}^{(i,j)}_{\mu_1,\mu_2,\nu_1,\nu_2,\vec{\alpha}}
= \{F_i=0\} \cap H_{i,\mu_1} \cap H_{i,\mu_2} \cap H_{j,\nu_1} \cap
H_{j,\nu_2} \cap
Z^{(i,j)}_{\vec{\alpha}},\]
where 
\[ Z^{(i,j)}_{\vec{\alpha}} = \bigcap_{\substack{k=1\\j\ne i, j}}^r \td_{k,\alpha_k}\]
is a toric subvariety of $X$ of codimension $r{-}2$.

\begin{figure}
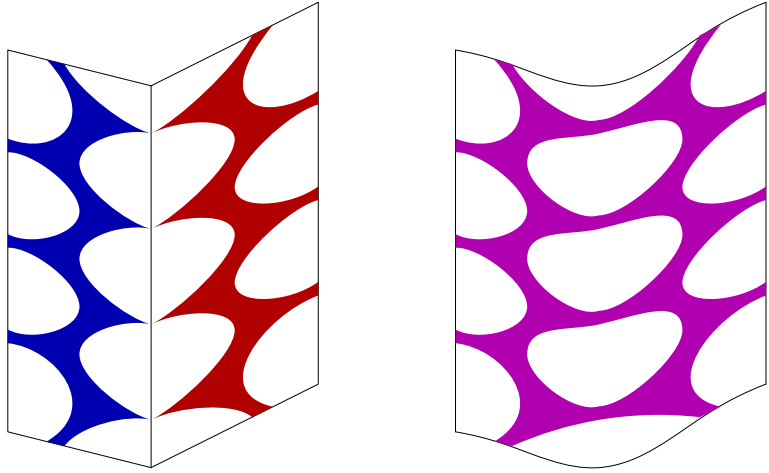

\begin{center}
\includegraphics[scale=0.5]{amoeb2.mps}
\qquad\qquad
\includegraphics[scale=0.5]{amoeb-junc.mps}
\caption{Components of  of $\pi(\widehat{C}^{(i)})$ meeting two-at-a-time
(left), and the deformation to a 
single component of $\pi(C^{(i)})$ (right).}
\label{fig:junction2}
\end{center}
\end{figure}

The reason for the extra varieties $\widehat{Q}$ is simple:  if we 
intersect $\widehat{C}^{(i)}_{\lambda_1,\lambda_2,\vec{\alpha}}$
with a toric divisor $H_a$, then we get something of the form
$\widehat{P}^{(i)}_{\lambda_1,\lambda_2,a,\vec{\alpha}}$ if $a$ belongs to the
$i^{\text{th}}$ part of the partition, but if $a$ belongs to the
$j^{\text{th}}$ part with $j\ne i$, then we get something of the form
$ \widehat{Q}^{(i,j)}_{\lambda_1,\lambda_2,\alpha_j,a,\vec{\alpha}}$.

\begin{figure}
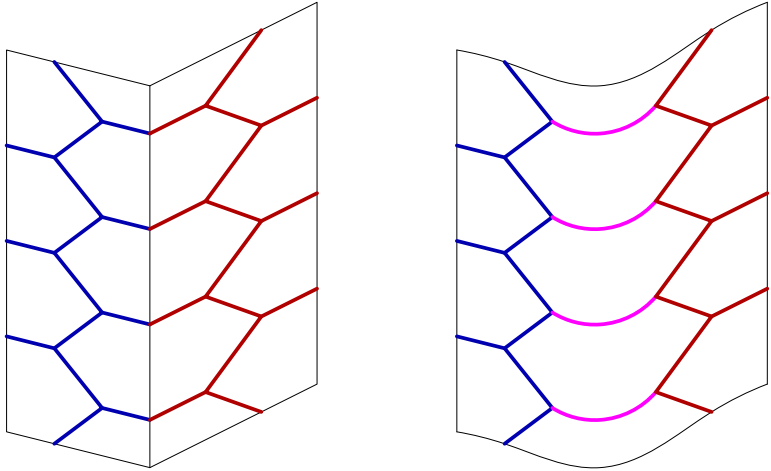

\begin{center}
\includegraphics[scale=0.5]{junction2.mps}
\qquad\qquad
\includegraphics[scale=0.5]{junctionm.mps}
\caption{Pieces  of $\widehat{D}^{(i)}$ meeting two-at-a-time
(left), and the discriminant underlying the
smoothing of a pair of components of
$\widehat{C}^{(i)}$ to a single component of $C^{(i)}$ (right).}
\label{fig:junction2-disc}
\end{center}
\end{figure}

Along intersections of type $\widehat{P}$, the components of 
$\pi(\widehat{C})$
are meeting three at a time (as is the case for the original singular 
locus $C$); this is illustrated in Figure~\ref{fig:junction3}.
On the other hand, along intersections of type $\widehat{Q}$,
the components of $\pi(\widehat{C})$ are only meeting two at a time,
as illustrated in the left side of Figure~\ref{fig:junction2}.
Passing back to the original $C$, components of $\pi(\widehat{C})$ will
merge along their two-at-a-time intersections, as illustrated
on the right side of  Figure~\ref{fig:junction2}.
The corresponding merger of combinatorial discriminants is illustrated
in Figure~\ref{fig:junction2-disc}.

Each component
$\widehat{C}^{(i)}_{\lambda_1,\lambda_2,\vec{\alpha}}$ is
a hypersurface in a toric variety 
$H_{i,\lambda_1} \cap H_{i,\lambda_2} \cap
Y^{(i)}_{\vec{\alpha}}$
of dimension $n-1$.
The moment map for that toric
variety maps 
$\widehat{C}^{(i)}_{\lambda_1,\lambda_2,\vec{\alpha}}$
to an open set 
$\pi(\widehat{C}^{(i)}_{\lambda_1,\lambda_2,\vec{\alpha}})$
in an $(n{-}1)$-dimensional polytope.
We expect that the image $\pi(C^{(i)}_{\lambda_1,\lambda_2})$
will degenerate to $\pi(\widehat{C}^{(i)}_{\lambda_1,\lambda_2})$,
which is a union of various open subsets of
$(n{-}1)$-dimensional polytopes, all contained within an
$n$-dimensional space (the moduli space of the $T^n$'s).
We formalize this as follows.

\bigskip

\begin{mainconjecturecontinued}
The image of each component of $C^{(i)}$ under $\pi$ is a deformation of 
a union of components of $\pi(\widehat{C}^{(i)})$, each of which is
an open subset of an $(n{-}1)$-dimensional polytope in $B$, with these
components meeting 
two-at-a-time.  The deformation smooths out the intersection points,
and is contained in an $(n{-}1)$-dimensional submanifold which is
a deformation of the union of the corresponding polytopes.

Upon choosing appropriate triangulation data of the dual polytope,
the image of the component of $C^{(i)}$ under $\pi$ retracts to a
codimension one subset of the $(n{-}1)$-dimensional submanifold
(a part of the combinatorial discriminant of $\pi$).  This
piece of the combinatorial discriminant is a deformation of
the union of
corresponding degenerate combinatorial discriminant pieces, each of which is
the retraction of a component 
of $\pi(\widehat{C}^{(i)})$ (within the corresponding 
$(n{-}1)$-dimensional polytope).

The degenerate combinatorial discriminant, as a codimension one
subset of the union of the polytopes associated to components of
$\widehat{C}$, coincides with the discriminant in 
the Haase--Zharkov construction 
\cite{HZ:T3III}.

\end{mainconjecturecontinued}

As stated in this part of the Main Conjecture,
if we begin with $\pi(\widehat{C}^{(i)}_{\lambda_1,\lambda_2})$,
it is a collection of amoebas in polytopes meeting two-at-a-time.
As the union of a pair of ambient polytopes smooths out, 
their amoebas will join smoothly, as illustrated 
in Figure~\ref{fig:junction2}.
At the same time, the corresponding varieties in $Z$ (which map to
those two components, and whose intersection is a double point
in complex codimension one) smooth to an irreducible variety.

One interesting observation is that the ambient toric varieties
of $\widehat{Q}^{(i,j)}_{\mu_1,\mu_2,\nu_1,\nu_2,\vec{\alpha}}$
and $\widehat{Q}^{(j,i)}_{\nu_1,\nu_2,\mu_1,\mu_2,\vec{\alpha}}$
are the same; in the first case, one also imposes $F_i=0$ while
in the second case one imposes $F_j=0$.  The varieties themselves
will only meet in codimension $4$ of $Z$, so when $n=3$ these
$\widehat{Q}$ sets are disjoint sets of points.  The fact that
they live in a common space leads to interesting possibilities
for linking and knotting of the combinatorial discrimimant, as
we shall point out explicitly in the conclusions.

\section{The combinatorial discriminant} \label{sec:CD}

Haase-Zharkov \cite{HZ:T3III} and Gross \cite{Gross:BB} produce a
combinatorial structure exhibiting some of the expected features of an
SYZ fibration.  From the combinatorial data determining a mirror pair
of complete intersection Calabi--Yau $n$-folds with trivial canonical
class in toric varieties, these authors construct a polytopal complex
$\Sigma$ with the topology of $S^n$, and two dual affine structures
defined on $\Sigma\setminus D$, where the combinatorial discriminant
$ D\subset \Sigma$
has real codimension $2$ (and is a trivalent graph when $n=3$).
These structures were studied in the hypersurface case in
\cite{susyT3}.  Here we review the construction for a complete
intersection, and show by explicit calculation in an example how the 
combinatorial discriminant $D$ is related to $\pi(C)$.  The complete 
intersection case introduces some additional subtleties into the
comparison.

Our example is the intersection $\P^5[4,2]$ of two hypersurfaces of
the indicated degrees in $\P^5$.  The Batyrev--Borisov construction
for this begins with the polytope associated to the ambient space.
This can be written in the lattice $N = \mathbb{Z}^5$ with basis $\{e_i\}_{i=1}^5$ as 
\begin{equation}\label{nablad}
\nabla^\vee = {\rm Conv}\left(0,e_0,e_1,e_2,e_3,e_4,e_5\right)
\end{equation}
where $e_0 = -\sum_{i=1}^5 e_i$.  
Our nef partition is $E^{(1)} = \{e_0,e_1,e_2,e_3\}$, $E^{(2)} =
\{ e_4,e_5 \}$, so that 
\begin{equation}\label{nsum}
\nabla^\vee = {\rm Conv}\left(\Delta^{(1)},\Delta^{(2)}\right)
\end{equation}
with 
\begin{equation}\label{deltas}
\Delta^{(1)} = {\rm Conv}\left(0,e_0,e_1,e_2,e_3\right)\qquad
\Delta^{(2)} = {\rm Conv}\left(0,e_4,e_5\right)\ .
\end{equation}
In this case, no triangulation of $\nabla^\vee$ is necessary.

The dual polytope $\nabla$ is a reflexive polytope in the dual lattice $M$.
This partition of $\nabla^\vee$ induces a decomposition of $\nabla$ as a Minkowski sum 
\begin{equation}
\nabla = \nabla^{(1)} + \cdots \nabla^{(r)}\ .
\end{equation}
This combinatorial construction is perhaps most intuitively understood
by considering $\nabla$ as the Newton polytope for a \CY\
hypersurface in $X$.  Explicitly, following the sign notation
of~\cite{HZ:T3III}, 
we associate to a point $m$ in the
dual lattice the monomial $M_m = \prod_a x_a^{1-\langle
  m,e_a\rangle}$ and thus lattice points in $\nabla$ are associated to
monomials with the same degree as $\prod_a x_a$.  
$\nabla^{(i)}$ then contains those points in $\nabla$ associated to
polynomials containing a factor of $\prod_{a\notin E_i} x_a$ and is
thus the Newton polytope for sections of the line bundle ${\cal
  L}_i$.  We have $\la\nabla^{(i)},\Delta^{(j)}\ra\le\delta^{ij}$.

In our example we find 
\begin{equation}\label{nablais}
\begin{aligned}
\nabla &= {\rm Conv}\{0,[1,1,1,1,1],[-5,1,1,1,1],[1,-5,1,1,1],\cr
&\quad [1,1,-5,1,1],[1,1,1,-5,1],[1,1,1,1,-5]\} \ ,\cr
\end{aligned}
\end{equation}
and the decomposition is 
$\nabla = \nabla^{(1)} + \nabla ^{(2)}$ with 
\begin{equation}\label{nablas}
\begin{aligned}
\nabla^{(1)} &= {\rm Conv}
\{0,[1,1,1,0,0],[-3,1,1,0,0],[1,-3,1,0,0],\cr
&\quad [1,1,-3,0,0],[1,1,1,-4,0],[1,1,1,0,-4]\}\cr
\nabla^{(2)} &= {\rm Conv}
\{0,[0,0,0,1,1],[-2,0,0,1,1],[0,-2,0,1,1],\cr
&\quad [0,0,0,-2,1,1],[0,0,0,-1,1],[0,0,0,1,-1]\}\cr\ .
\end{aligned}
\end{equation}
Our example enjoys a large symmetry, and this is manifest in the
decomposition above.  
Upon inspection one sees that $\nabla^{(1)}$ is
obtained from 
$\nabla^{(2)}$ by translation (by $[1,1,1,-1,-1]$) and
doubling.  This reflects the fact that we are in $\mathbb{CP}^n$.  
We label the
vertices of $\nabla^{(i)}$ as $v^{(i)}_a,\ a=0,\ldots, 5$, reflecting
this symmetry.

Mirror symmetry in this context is a combinatorial duality
\cite{borisov-ms,MR1416334,Batyrev-Borisov-1994}.  A triangulation of
the reflexive polytope 
\begin{equation}
\Delta^\vee = {\rm Conv}(\nabla^{(1)},\ldots,\nabla^{(r)})
\end{equation}
determines a dual toric variety, nef partition, and a mirror family of
complete intersection \CY\ spaces, and
\begin{equation}
\Delta = \Delta^{(1)} + \cdots \Delta^{(r)}\ .
\end{equation}

In our example, setting $\Delta^\vee = {\rm Conv}\left(\nabla^{(1)},\nabla^{(2)}\right)$ we find 
\begin{equation}\label{deltais}
\Delta = \Delta^{(1)} + \Delta^{(2)} 
= {\rm Conv}\{0,e_a,e_I+e_\alpha\}\quad I\in \{0\ldots 3\};\ \alpha\in\{4,5\}\ .
\end{equation}

We are interested in studying $Z$ deep inside its K\" ahler cone
 -- determined by a triangulation of $\nabla^\vee$ -- and near a complex
structure limit point with maximal unipotent monodromy.  This
corresponds via the monomial-divisor mirror map to a point deep in the
K\"ahler cone of the mirror -- determined by a triangulation of
$\Delta^\vee$.   We thus consider both polytopes to be triangulated.

The construction of \cite{HZ:T3III} produces a polytopal complex $\Sigma$ in
$\Delta\times\nabla$ with the topology of $S^3$.  $\Sigma$ is a subdivision
of $|\Sigma| = \{(n,m)\in\Delta\times\nabla : \la n,m\ra = r\}$
determined by the triangulations of the respective polytopes.  We can
see some of the structure before triangulating.
For a proper face $t\subset \Delta^\vee$ we define $t^{(i)} =t\cap\nabla^{(i)}$.  
A face will be called {\sl transversal\/} if none of these are empty.
Minimal transversal faces are $r-1$ simplices with one
vertex in each $\nabla^{(i)}$.
For a transversal $t$, $t_\nabla = t^{(1)} + \cdots +
t^{(r)}\subset\nabla$ is a proper face of
$\nabla$ \cite{Gross:BB}.   For a minimal face this is a vertex of $\nabla$.
The dual face to this $s\in\nabla^\vee$ is
transversal, i.e.~$s^{(i)} =s\cap\Delta^{(i)}$ are all non-empty, and 
$s_\Delta = s^{(1)} +\cdots + s^{(r)}$ is the face of $\Delta$  dual to $t$.  
Faces $s\subset\Delta^\vee$ and $t\subset\nabla^\vee$ so related are
{\sl adjoint\/} and we have ${\rm dim} t_\nabla + {\rm dim} s_\Delta =
n$.  Clearly here $\la s_\Delta,t_\nabla\ra = r$ and 
$|\Sigma|$ is the collection $s_\Delta\times t_\nabla$ for adjoint
$(s,t)$.  

In our example, 
a face $s$ of the simplex $\nabla^\vee$ (of any
dimension) is transversal if it contains one vertex from each of the
$\Delta^{(i)}$.  For example, there are 8 irreducible transversal faces
$s^{I\alpha} = {\rm Conv}\{e_I,e_\alpha\}$, and these produce
eight of the vertices of $\Delta$: $s^{I\alpha}_\Delta = e_I +
e_\alpha$.
With a hopefully obvious notation the transversal faces of
$\nabla^\vee$ 
and the faces of $\Delta$ 
that they produce are 
\begin{equation}\label{nablavfaces}
\begin{aligned}
s^{IJ\alpha}_\Delta &= {\rm
Conv}\{s^{I\alpha}_\Delta,s^{J\alpha}_\Delta\}\cr
s^{I45}_\Delta &= {\rm Conv}
\{s^{I4}_\Delta,s^{I5}_\Delta \}\cr
s^{IJK\alpha}_\Delta &= {\rm Conv}\{s^{IJ\alpha}_\Delta,
s^{IK\alpha}_\Delta,s^{JK\alpha}_\Delta\}\cr
s^{IJ45} &= {\rm Conv}\{s^{I4}_\Delta,s^{I5}_\Delta,
s^{J4}_\Delta,s^{J5}_\Delta\}\cr
s^{0123\alpha} &= {\rm Conv}\{s_\Delta^{0\alpha},s_\Delta^{1\alpha},s_\Delta^{2\alpha},s_\Delta^{3\alpha}\}\ .\cr
\end{aligned}
\end{equation}
Note that $\Delta$ is not a simplex.  
The two-dimensional faces we find are either triangles ($\dim s^{(1)}=2,\
\dim s^{(2)}=0$) or quadrilaterals
($\dim s^{(1)}=\dim s^{(2)}=1$).

The six irreducible transversal faces of $\Delta^\vee$
are $t^{aa} = {\rm Conv}\{v^{(1)}_a,v^{(2)}_a\}$, and they lead to the
six vertices of $\nabla$
$t^a_\nabla = v^{(1)}_a + v^{(2)}_a$.  The faces of dimension one and
two are
\begin{equation}\label{deltavfaces}
\begin{aligned}
t^{IJ}_\nabla &= {\rm Conv}\{t^I_\nabla,t^J_\nabla\}\cr
t^{I\alpha}_\nabla &= {\rm Conv}\{t^I_\nabla,t^\alpha_\nabla\}\cr
t^{IJK}_\nabla &= {\rm Conv}\{t^I_\nabla,t^J_\nabla,
t^K_\nabla\}\cr
t^{IJ\alpha}_\nabla &= {\rm Conv}\{t^I_\nabla,t^J_\nabla,
t^\alpha_\nabla\}\cr
t_\nabla^{IJK\alpha} &= {\rm Conv}\{t^I_\nabla,t^J_\nabla,t^K_\alpha,
t^\alpha_\nabla\}\cr\ .\cr
\end{aligned}
\end{equation}
Sixteen of the twenty triangular two-faces of $\nabla$ are transverse.

The large symmetry of our example makes the duality straightforward:
$s^{\cal I}_\Delta$ and $t^{\cal J}_\nabla$ are adjoint faces
precisely when $\{0,\ldots,5\} = {\cal I}\amalg{\cal J}$.  More
generally, if $s\subset s'$ is a transversal face of $s'$ (implying
$s'$ is transversal) and if $t$ is adjoint to $s$ then the adjoint
$t'$ to $s'$ is a face of $t$.

The combinatorial structure is now manifest.
The four faces $t_\nabla^{IJK}$, each adjoint to an edge
$s_\Delta^{L45}$, form a tetrahedron.  The twelve remaining 
faces $t_\nabla^{IJ\alpha}$, adjoint to edges $s_\Delta^{KL\beta}$,
 fall into two groups of six (by $\alpha$). 
Along each of the six edges $t_\nabla^{IJ}$, two of faces of the first type
meet two of the second type, one from each family.  In accordance with
the expected duality, these edges are adjoint to the quadrilateral
faces $s_\Delta^{KL\alpha\beta}$.  Each edge of $s$ is adjoint to one of the
intersecting faces.
Along the remaining eight edges $t_\nabla^{I\alpha}$ the members of each
family meet three at a time; these edges are adjoint to the triangular
faces $s_\Delta^{JKL}$.

A triangulation of $\nabla^\vee$ and $\Delta^\vee$ determines, as
mentioned above, a suitable limiting point.  This induces
subdivisions
$S,T$ of the boundaries 
$\partial\nabla^\vee$ and $\partial\Delta^\vee$ respectively.
This in turn induces a subdivision of $\Sigma$ into cells of the form
$(\sigma_\Delta,\tau_\nabla)$ where $\sigma\subset s$ and $\tau\subset
t$ are contained in pairs of adjoint faces.
The construction of \cite{HZ:T3III,Gross:BB}
produces a trivalent graph in $\Sigma$ representing the
discriminant of the fibration in the large radius and large
complex structure limit.

To construct the discriminant we take the barycentric subdivision
associated to $\Sigma$.  An adjoint pair $(\sigma,\tau)$ is {\sl
  smooth\/} if $\dim\sigma^{(i)}\cdot\dim\tau^{(i)} = 0,\ \forall i$.
Vertices of the discriminant are associated to  {\it non-smooth\/}
adjoint pairs.  
Clearly, if $\sigma$ or $\tau$ is a minimal transversal face the pair
is smooth.  Non-smooth cells will thus be associated to adjoint pairs
for which $\dim\sigma_\Delta$ and $\dim\tau_\nabla$ are both positive.  
As a result, both are contained within faces of $\Delta$,
resp. $\nabla$, of dimension at most $n-1$.

For the purpose of constructing the
discriminant, we can thus limit ourselves to the triangulation of
these faces.  A transversal $r$-simplex $s$ in the
boundary of $\nabla^\vee$ producing an edge of $\Delta$
will have two vertices $e_{i,\lambda_1}, e_{i,\lambda_2}$ lying in one
component $\Delta^{(i)}$ and one vertex $e_{j,\nu_j}$ in each of the other
components.  In $X$ this corresponds to the codimension-$r+1$ toric
subvariety $x_{i\lambda_1}=x_{i,\lambda_2} = x_{j,\nu_j} = 0$.   
In the terminology of section~\ref{sec:NP} this is the toric
subvariety $H_{i,\lambda_1}\cap H_{i,\lambda_2}\cup Y^{(i)}_{\vec\alpha}$
with $\alpha_j = \nu_j\ \forall j\neq i$.

The dual face $t$ of $\Delta^\vee$ will be
adjoint to $s$ and will produce an $n-1$ dimensional face of $\nabla$.
Faces of $t$ will be adjoint faces containing $s$
as a face, so all non-smooth cells will be contained in the
triangulation of these $t$.  The vertices of $D$ will correspond to
the non-smooth cells in the triangulation.

The toric subvariety $H_{i,\lambda_1}\cap H_{i,\lambda_2}\cup Y^{(i)}_{\vec\alpha}$
maps to the associated face of $\nabla^\vee$ under the moment map
$\mu:X\to N_{\R}$, and the combinatorial version of the conjecture is
that the projection of the component
$C^{(i)}_{\lambda_1,\lambda_2,\vec\alpha}$ retracts to the graph so
constructed.   We will not prove this in general.

In our example, the fact that $\nabla^\vee$ requires no subdivision
makes things simple.  The transversal two-faces producing edges of
$\Delta$ are dual to the three-faces of $\Delta^\vee$ producing the
two types of two-dimensional faces of $\nabla$ described above.  Our
job is to triangulate these.   The existence of a subdivision of the
polytope consistent with our choices is supported by the consistency
of our results, although we have not proved it.  

\begin{figure}
\begin{center}
\includegraphics[scale=0.75]{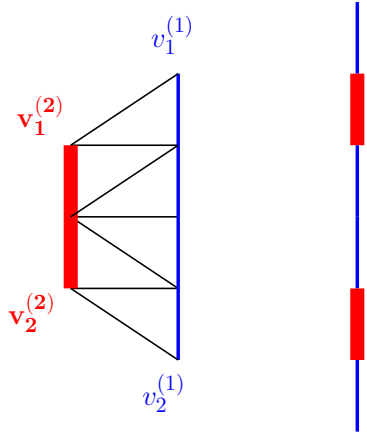}
\end{center}
\caption{Triangulation of a trapezium leading to the subdivision of an
interval.}\label{fig:XX}
\end{figure}

Consider first the one-dimensional cell $\tau^{01}_\nabla$.  The
two-dimensional transverse face producing this is the trapezoid with
vertices $v^{(1)}_0,v^{(1)}_1,v^{(2)}_0,v^{(2)}_1$.  Subdividing this
as in Figure~\ref{fig:XX} we find six cells and the resulting division of
$\tau^{(01)}$ is indicated.  The six cells into which the segment is
subdivided are colored to indicate that four of them have
$\dim\tau^{(1)}=1,\dim\tau^{(2)}=0$ while two have 
$\dim\tau^{(1)}=0,\dim\tau^{(2)}=1$.  We assume that we can make this
subdivision for each of the $14$ transversal edges.

\begin{figure}
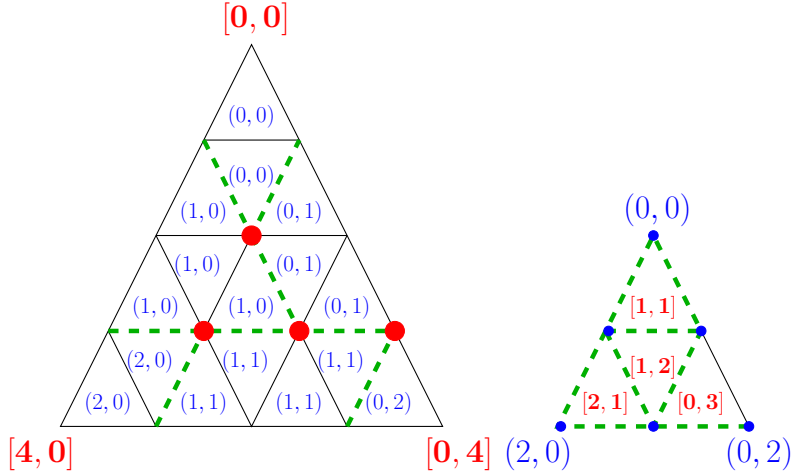

\begin{center}
\includegraphics[scale=0.5]{4triangle5Bb.mps}
\includegraphics[scale=0.5]{2triangle3aB.mps}
\end{center}
\caption{How to match the triangles of different sizes.  The large
(red) vertices on the left are paired with correspondingly-labeled
faces on the right, and the small (blue) vertices on the right are
paired with faces on the left.
The dashed (green) edges pair up between the triangles, to span simplices
between them.  (For example, the edge between $[2,1]$ and $[1,2]$
pairs with the edge between $(1,0)$ and $(1,1)$.)
}\label{fig:triangulation}
\end{figure}

\begin{figure}[ht]
\begin{center}
\includegraphics[scale=0.75]{prismaBB.mps}
\caption{Triangulating a prism.  The upper triangle is a sub-triangle
of the right side of Figure~\ref{fig:triangulation}, while the
lower triangle is a sub-triangle of the left side of Figure~\ref{fig:triangulation}.}\label{fig:prism}
\end{center}
\end{figure}

\begin{figure}
\begin{center}
\includegraphics[scale=0.5]{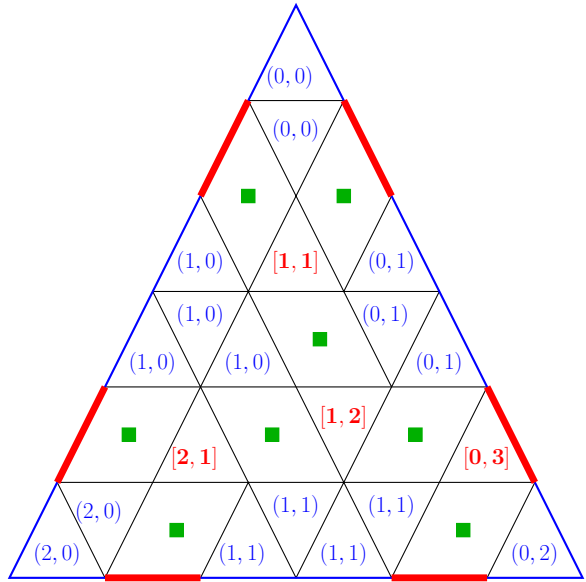}
\end{center}
\caption{The corresponding large triangle, with both types of faces,
labeled in bold (red) type and non-bold (blue) type.
The square (green) vertices mark faces which are the Minkowski sums of
a pair of dashed (green) edges from the previous figures.
}\label{fig:largetriangle}
\end{figure}

We now need to consider the two-dimensional cells of either type.  
These will all be
equilateral triangles with sides of length six as we have seen above.
Each is obtained from a three-dimensional face of $\Delta^\vee$
intersecting $\Delta^{(i)}$ in 
two parallel triangles of side lengths four and two, as in
Figure~\ref{fig:triangulation}.  The region between those triangles
consists of triangular prisms which must be further subdivided.
Our subdivision includes three kinds of simplices.  
The first type, indicated
in boldface (red) type in the figure, meets the larger triangle at a vertex
and the smaller in a face and so has dimensions
$(\dim\tau^{(1)},\dim\tau^{(2)})=(0,2)$.  
The second type, indicated with ordinary (blue) type in the figure,
inverts this and has dimensions $(2,0)$.  A third type,
indicated by the dashed (green) sides, meets each triangle along a side so has
dimensions $(1,1)$.  Figure~\ref{fig:prism} 
shows how these three types arise in subdividing a prism.
As with the interval in Figure~\ref{fig:XX}, 
we obtain an induced subdivision of
the two-dimensional face by taking Minkowski sums.  The sum of a
triangle and a point leads to triangles of the two types shown in
Figure~\ref{fig:largetriangle}. 
The sum of two intervals is a parallelogram, leading to the
(green) cells that pair up two dashed edges.  
The figure illustrates that 
this is compatible with our
subdivision of the edges, in that the two segments of each side
that are associated to $\nabla^{(2)}$ are in the positions in which we
need them.

\begin{figure}
\begin{center}
\includegraphics[scale=0.5]{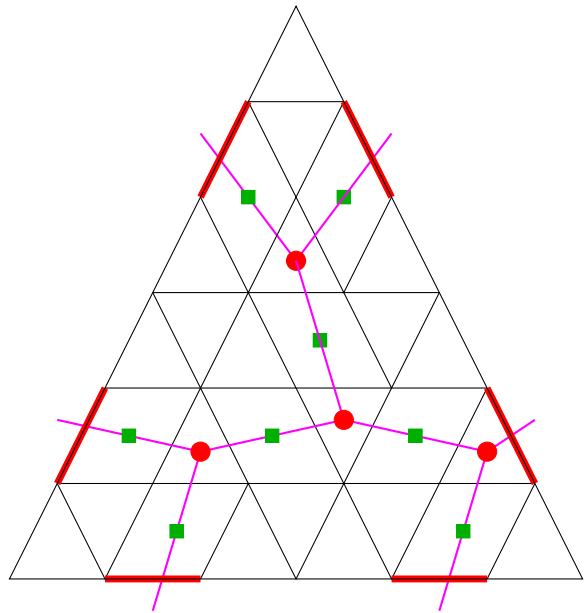}
\end{center}
\caption{The triangulation with large (red) and square (green) vertices.
}\label{fig:redgreen}
\end{figure}

\begin{figure}
\begin{center}
\includegraphics[scale=0.5]{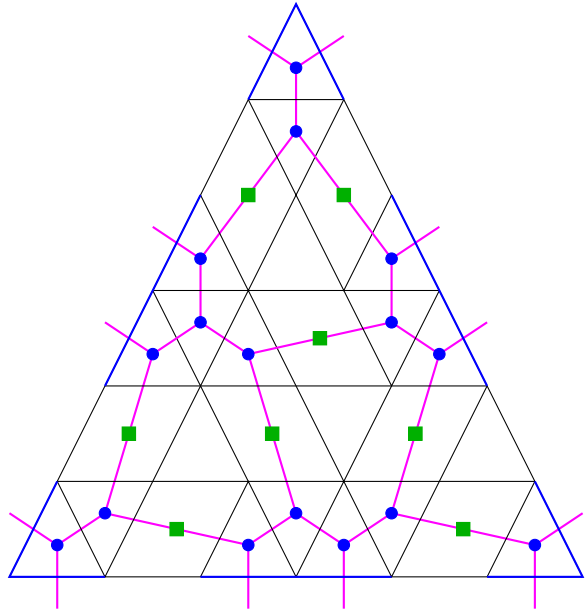}
\end{center}
\caption{The triangulation with small (blue) and square (green) vertices.
}\label{fig:bluegreen}
\end{figure}

With our subdivision in hand we can now find the non-smooth pairs
whose centers will form the vertices of our graph for the
discriminant.  Non-smooth cells with
$(\dim\sigma_\Delta,\dim\tau_\nabla) = (1,2)$ occur on the large
triangles just discussed.  As mentioned above, there are two types of
these faces.  Faces of the first type, $\tau_\nabla^{IJK}$, are paired with
one-dimensional cells such as $\sigma_\Delta^{L45}$.  These have
dimensions $(0,1)$ so non-smooth pairs will result when the subdivided
$\tau$ has nonzero dimensional intersection with $\nabla^{(2)}$,
i.e., for the bold (red) and dashed (green) cells 
as indicated in Figure~\ref{fig:redgreen}.  

These faces correspond, as discussed above, to the four toric subvarieties
$H_{1,L}\cap H_{2,4}\cap H_{2,5}$.  In the limit in which $F_1 =
x_0x_1x_2x_3$ the (single) component $C^{(2)}_{45}$ restricts to the
union of these as a plane conic in each $\P^2$.    
Note that here $\widehat P^{(2)}$ is empty.  

The twelve faces of the other type $\tau_\nabla^{IJ\alpha}$, are
paired with one-dimensional cells $\sigma_\Delta^{KL\beta}$ with
dimensions $(1,0)$.  Here non-smooth pairs arise from cells with nonzero
dimensional intersection with $\nabla^{(1)}$, i.e., the non-bold (blue) and 
dashed (green)
cells as indicated in Figure~\ref{fig:bluegreen}.   

These correspond to the twelve toric subvarieties $H_{1,K}\cap
H_{1,L}\cap H_{2,\beta}$.  In the limit in which $F_2 = x_4 x_5$ the
components of $C^{(1)}$ restrict to the union of these as plane
quartics in each $\P^2$.

Note that
each of the quadrilateral (green) cells adjoins two non-smooth cells
(along a non-smooth edge) and two smooth cells (along a
smooth edge) so that these cells lead to bivalent vertices.  The
adjacent smooth cells show that the monodromy about the two
edges meeting at such a vertex is the same.

We can now compute the Euler character of $Z$, which is predicted to
be the difference of the number of positive and negative vertices of
$D$.  Each of the four faces $\tau^{IJK}$ has four negative vertices
(note that the bivalent vertices can be smoothed and do not contribute
to the Euler character),  while each of the twelve faces
$\tau^{IJ\alpha}$ has 16 negative vertices, for a total of 208 negative
  vertices.  

The positive vertices are associated to non-smooth
cells with $(\dim\sigma_\Delta,\dim\tau_\nabla) = (2,1)$.  These arise,
as discussed above,
along the edges of the large faces. (non-smooth internal edges connect
two non-smooth cells so the associated vertex is a simple double
point and does not contribute to our Euler characteristic).  As
mentioned above, there are two types of edges.  At an edge 
$\tau_\nabla^{IJ}$ two faces of each type meet.
Since the associated $\sigma_\Delta^{KL45}$ has dimensions $(1,1)$ all cells
are non-smooth.  The $\sigma_\Delta$ cells are quadrilaterals and as was the
case above, each connects two non-smooth cells (belonging to faces of
identical type) along a non-smooth edge and two smooth cells
(belonging to the faces of the other type) along a smooth edge.  Thus
the six quadrilaterals along each of these large edges are all
bivalent vertices (oriented in two different ways) connecting faces of
the same type.  All of this contributes nothing to the Euler
characteristic.

These edges are associated to the toric subvarieties $H_{1,K}\cap
H_{1,L}\cap H_{2,4}\cap H_{2,5}$.  As discussed above, we see that
these contain both the intersections $\widehat Q^{(1,2)}_{KL45}$ and $\widehat
Q^{(2,1)}_{45KL}$ (for $r=2$ the $Z^{(1,2)}$ are just $X$).  
These are generically disjoint two-at-a-time
intersections of pairs of components of $C^{(1)}$ and $C^{(2)}$.

The other type of edge $\tau_\nabla^{I\alpha}$, 
along which three faces of the second type
meet, is paired with a cell 
$\sigma^{JKL\beta}$ of dimensions $(1,0)$.  Segments of the edge
with nonzero dimensional intersection with $\nabla^{(1)}$ (the four
blue segments) are
non-smooth.  These are precisely the edges bounding non-smooth cells,
and so each such edge contributes four trivalent vertices. There are
eight such edges and hence 32 positive vertices.

These are associated to  $H_{1,J}\cap H_{1,K}\cap H_{1,L}\cap H_{2,\beta}$
and contain the intersections $\widehat P^{(1)}_{JKL\beta}$.

Adding up the contributions to the Euler characteristic we find 
$32 - 208 = -176$ as expected.

\section{Concluding Remarks}

Strominger, Yau, and Zaslow applied physics arguments to 
a string-theoretic version of 
 mirror symmetry, obtaining a
remarkable prediction about mathematics.
An understanding of the structure of the fibration whose existence
they predicted
could provide insights into both \CY\ geometry and mirror symmetry.
The geometric conjectures of \cite{susyT3} are geometrically natural
but some of the strongest evidence for them was their compatibility
with the combinatorial constructions of 
\cite{HZ:T3I,HZ:T3II} (which provided a general formulation of the 
original constructions
given by Gross and by Ruan).
The
extension of these combinatorial
constructions to complete intersections is quite
nontrivial.  The fact that the conjecture extends naturally in a way
that appears compatible with this more elaborate construction is
encouraging.  

One new feature of the complete intersection construction is that in
general the components of the combinatorial discriminant, graphs in
$S^3$ for $n=3$, will be linked and (possibly) knotted.  We have not
attempted to compute this in general.

Note that in this paper we have not proved the compatibility of our 
conjectured construction with the complete intersection combinatorics
but have simply checked it in an example.  It should be possible
to show that for any choice of triangulation data that fully specifies
the discriminant in the combinatorial construction of \cite{HZ:T3III},
the images under $\pi$ of the components of $\widehat{C}$ do indeed
retract to pieces of the discriminant; we have not attempted this.
In particular, the detailed structure of $D$
depends on a choice of triangulation, or equivalently a choice of a
limiting point of maximally unipotent monodromy in moduli space. 
Different choices of triangulation of $\nabla^\vee$ are related by
flop transitions, and understanding
the way the fibration changes in these would be useful.

Mirror symmetry was an essential ingredient in the original argument
for the existence of a fibration and one of the most intriguing
predictions is that the fibrations of two members of a mirror pair are
related by a simple duality.  One attractive feature of the
combinatorial construction is that it implements this duality very
naturally.  The combinatorial construction for the mirror of $Z$ is
obtained by exchanging the roles of $\Delta$ and $\nabla$.   The
construction of $\Sigma$ is clearly invariant under this.  Moreover
the construction leads naturally to two dual integral affine
structures sharing the same discriminant. 

The way in which the conjecture is compatible with this symmetry is
obscured by the simple nature of our example, so it perhaps bears
mention.  In our example, $\nabla^\vee$ was a simplex and the
combinatorial structure was essentialy determined by the triangulation
of $\Delta^\vee$:  all the cells in a triangulation of a face of
$\Delta^\vee$ were paired with the same face of $\nabla^\vee$.   In a
more general example, the face of $\nabla^\vee$ would also be
triangulated.  For example, in the combinatorial construction for the 
mirror of our example -- a codimension two complete intersection in 
a toric variety determined by a triangulation of $\Delta^\vee$ --  if we
choose the triangulation used here, the red triangles in faces would
be associated to the ambient space of $\widehat P^{(2)}$, 
the blue triangles to the ambient space of $\widehat
P^{(1)}$, and the quadrilateral cells to the common ambient space of
$\widehat
Q^{(1,2)}$ and $ \widehat Q^{(2,1)}$.  

\medskip
\noindent{\bf Acknowledgements}:  We are grateful to P.S~ Aspinwall,
M.~Bertolini, R.~Casta\~no-Bernard, C.~Haase, and X.~de la Ossa for 
helpful conversations,
and to D.-E.~Diaconescu for collaboration in the early stages of this
work.  MRP thanks the Mathematics department at UCSB for gracious
hospitality during essential phases of this work.  
DRM is supported by NSF grant PHY-1307513, and
MRP is supported by 
NSF grant PHY-1217109.
Any opinions, findings, and conclusions or
recommendations expressed in this material are those of the authors
and do not necessarily reflect the views of the National Science
Foundation.

\providecommand{\bysame}{\leavevmode\hbox to3em{\hrulefill}\thinspace}
\providecommand{\MR}{\relax\ifhmode\unskip\space\fi MR }
\providecommand{\MRhref}[2]{%
  \href{http://www.ams.org/mathscinet-getitem?mr=#1}{#2}
}
\providecommand{\href}[2]{#2}

\end{document}